\documentclass{tran-l}
\usepackage{graphicx}
\usepackage{latexsym}
\usepackage{amssymb}
\def\bbbr{{\rm I\!R}} %reelle Zahlen
\def\bbbc{{\mathchoice {\setbox0=\hbox{$\displaystyle\rm C$}\hbox{\hbox
to0pt{\kern0.4\wd0\vrule height0.9\ht0\hss}\box0}}
{\setbox0=\hbox{$\textstyle\rm C$}\hbox{\hbox
to0pt{\kern0.4\wd0\vrule height0.9\ht0\hss}\box0}}
{\setbox0=\hbox{$\scriptstyle\rm C$}\hbox{\hbox
to0pt{\kern0.4\wd0\vrule height0.9\ht0\hss}\box0}}
{\setbox0=\hbox{$\scriptscriptstyle\rm C$}\hbox{\hbox
to0pt{\kern0.4\wd0\vrule height0.9\ht0\hss}\box0}}}}
\def\bbbz{{\mathchoice {\hbox{$\mathsf\textstyle Z\kern-0.4em Z$}}
{\hbox{$\mathsf\textstyle Z\kern-0.4em Z$}}
{\hbox{$\mathsf\scriptstyle Z\kern-0.3em Z$}}
{\hbox{$\mathsf\scriptscriptstyle Z\kern-0.2em Z$}}}}

\theoremstyle{definition}

\theoremstyle{remark}

\numberwithin{equation}{subsection}
\begin{document}

\title[Towards a Broader View of Theory of Computing]
 {Towards a Broader View of Theory of Computing - Part 1 }

\author{Narendra Karmarkar}

\address{Department of Electrical Engineering, Indian Institute of Technology, Bombay, India}

\email{narendrakarmarkar@yahoo.com}

%\subjclass{}

%\keywords{ }

\date{December 8, 2014 }

%\dedicatory{}

%\commby{Daniel}

% -----------------------------------------------------------------------------

\begin{abstract}

Beginning with the projectively invariant method for linear programming, interior point methods have led to powerful algorithms for many difficult computing problems, in combinatorial optimization, logic, number theory and non-convex optimization. Algorithms for convex optimization benefitted from many pre-established ideas from classical mathematics, but non-convex problems require new concepts. \\ \\
Lecture series I am presenting at the conference on Foundations of Computational Mathematics, 2014, outlines some of these concepts-- computational models based on the concept of the continuum, algorithms invariant w.r.t. projective , bi-rational , and bi-holomorphic  transformations  on co-ordinate representation, extended proof systems for more efficient certificates of optimality, extensions of Grassmann’s extension theory, efficient evaluation methods for  the effect of exponential number of constraints , theory of connected sets  based on graded connectivity, theory of curved spaces adapted to the problem data,  and concept of “ relatively” algebraic  sets in curved space.\\ \\
Since this conference does not have a proceedings, the purpose of this article is to provide the material being presented at the conference in more widely accessible form.
\end{abstract}

% -----------------------------------------------------------------------------
\maketitle
% -----------------------------------------------------------------------------

\section{Introduction}

 In Part 1 of this lecture series we discuss two topics : 
 \begin{itemize}
   \item Computational models based on the concept of the continuum
   \item Extended proof systems for more efficient certificates of optimality
 \end{itemize}
  In the part 2 of this lecture series, we describe efficient evaluation methods for  the effect of exponential number of constraints in the context of continuum based algorithms.We illustrate these ideas by considering concrete examples of finding maximum independent set in a graph and the satisfiability problem. Objective function for these problems is treated by non-convex optimization methods covered in part 3.

\section{\textbf{Models of Computation}}
\subsection{\textbf{Introduction }}
A model of computation provides mathematical abstractions of basic data objects and operations on those objects, available as building blocks. This effectively decouples design of algorithms for complex tasks from lower level details of how the underlying hardware implements the basic primitives. The Turing machine model uses strings of 0’s and 1’s and finite state machines. Careful study of the work of early pioneers – Turing,Von Neumann and Godel – shows that they were acutely aware of the limitations of this model for comprehensive understanding of the fundamental aspects of computation.  BSS model(Blum, Shub, Smale) uses real or complex numbers as data objects and algebraic operations ( including comparisons in real case ). This is more natural for many algorithms in numerical analysis, whereas Turing Machine model might seem more appropriate for discrete applications.

\subsection{\textbf{Applications of computing}}
Various applications of computing come in both flavours- discrete and continuous. e.g. Many business and information technology applications of computing clearly use discrete models.

At the same time, there are many areas of scientific and engineering applications based on continuous models. A few such examples are :
\begin{itemize}
\item Numerical simulation of natural phenomena or engineered systems based on differential equations
\item Estimation of probabilities, assessing strength of association between various entities in social networks
\item Interior point algorithms in optimization which are based on embedding discrete configurations in a multidimensional continuous space and constructing trajectories converging to the solution
\end{itemize}
Ideas presented in this paper are primarily motivated by such applications.
A model of computation is used to map applications of interest on physical machines \\

\subsection{ \textbf{Examples of physical computers}}
\begin{itemize}
\item Standard digital computers
\item Analog computers explored in the past e.g. Shannon's differential analyzer
\item Quantum computers
\item Biological systems that appear to do information processing / computing
\item Computers inspired by biological systems and partially mimicking some aspects of those systems \\
 e.g. neural networks, neuromorphic computing, etc.

\end{itemize}

\subsection{ \textbf{The need for a broader approach -- Early Pioneers' views } }
\begin{itemize}
\item \textbf{Alan Turing} was interested in understanding the origin of intelligence in biological systems
He believed it is due to some form of computing happening in these systems.
However, there doesn't seem to be anything similar to Turing machines in these systems.
Many years after developing the discrete model of computation,
Turing started exploration of partial differential equations modelling biological functions
These are clearly continuum based models of biological systems.

\item \textbf{Von Neumann} defined architecture of a digital computer using many ideas from Turing Machine model.
However, he was particularly critical of the limitations imposed on the
theory of automata by its foundations in formal logic, combinatorics.
He articulated the need for a detailed, highly mathematical, and more
specifically, analytical theory of computation based on mathematical analysis.
 \item \textbf{Shannon} worked on analog computer called differential analyzer. However, direct implementation of differentiation in analog system is highly error prone as it involves subtraction of two nearly equal continuous quantities.
\end{itemize}
\newpage
\subsection{\textbf{ Physical devices that process information in continuum form}}

\begin{itemize}
\item \textbf{Role of Integration and Integral Transform}

\begin{itemize}
\item It's well known that for certain non-linear differential equations which are difficult to solve numerically, it helps to convert them into equivalent integral equations.
\item Biological systems seem to use integral transform. e.g. Human ear essentially computes the magnitude of Fourier Transform of speech signal.The transform of the derivative is then obtained by scaling, a simpler algebraic operation than differentiation in analog setting.
\item Some optical computers are implementing Fourier Transform directly for computing differential operators required in computational fluid dynamic
\item Form of integration or convolution restricted to non-negative real functions is easier to provide in a physical system.  
\end{itemize}

\item \textbf{Robustness }

\begin{itemize}
\item In digial systems, the mapping from \\
  \\
 \textbf{Physical states} / \textbf{quantities} $\hspace{7.0pt}\longrightarrow \hspace{7.0pt}$ \textbf{information states}/\textbf{quantities}\\
 \\
 is many to one, In fact, each information state has infinite pre images.

\item This is to gain robustness in the face of small variations in physical states, due to noise, thermal effects, variations in manufacturing etc.
\item However, this does not necessarily require information states to be discrete. e.g. An ideal low pass filter also provides infinity to one mapping but the output signal still belongs to the continuum.
\end{itemize}

\end{itemize}

It is important to include these two mechanisms when constructing machines
that support continuum computing

\subsection{\textbf{Intertwining of discrete and continuous models}}
Logical and physical processes underlying computation involves several levels of abstractions, both continuous and discrete. Let us examine a specific application -- Solution of a non-linear differential equation.

\begin{itemize}
  \item This application and algorithm used for solution are both based on continuous model.
  \item  Real numbers in the continuous model are approximated by floating point arithmetic for mapping onto the discrete Turing machine model.
  \item  The computer implementing discrete Turing Machine model is actually a network of transistors. As a circuit, it is modeled by a system of non-linear differential equations. Restricted interpretation is applied to continuous quantities. e.g. Continuous time is divided into clock cycles bigger than settling time of the transient solution to the differential equations, to create discrete time steps.Voltage greater than certain threshold voltage represents logical one. In this way, the discrete logical model is supported by the underlying circuit-level continuous model.
  \item  This continuous circuit model is itself approximating underlying discrete phenomenon :
  \begin{itemize}
    \item Continuous, fluid-like model of current flow, for approximating the statistical mechanics of large number of discrete particles (electrons and holes in a semiconductor)
    \item Continuous energy bands( valence and conduction bands) representing large number of closely spaced quantized energy levels
  \end{itemize}
  \item  The underlying quantum mechanical phenomenon is modelled mathematically by the concept of Hilbert space which is based on the continuum concept. Similarly, more recent and on-going model of physics is based on string theory, which again requires the concept of the continuum.
\end{itemize}

%\begin{figure}
%%\includegraphics[scale=0.4]{}
%\end{figure}
%\begin{figure}
% \includegraphics[scale=0.4]{table1}
%\end{figure}

\subsection{ \textbf{Towards a broader view of theory of computing}}

\paragraph{ Given this intertwining of discrete and continuous models from top to bottom, it would be more illuminating to take a broader view of theory of computing. In addition to the discrete Turing machine model, one also needs to explore continuous models for a comprehensive understanding of the fundamental aspects of computation.  BSS model uses real or complex numbers as data objects and algebraic operations ( including comparisons in real case ). This is more natural for many algorithms in numerical analysis. }
\paragraph{ Various computing models can be organized in a similar way as Cantor had organized infinite sets -- by cardinal number of the set of all possible machines and data objects in the model. Staying within the same cardinal number, a more powerful approach is to use further extension, e.g real analytic functions or algebraic closure of meromorphic functions over suitable domains. Operations include algebraic as well as analytic operations. i.e. integration and differentiation regarded as binary operations. ( specification of the contour of integration is one of the input operands) }
\paragraph{ All such models that use both algebraic and analytical operations on continuous data objects are collectively referred to as continuum computing.Our research program is aimed at exploration of what a continuum view of computing might suggest for science of algorithms, relative difficulty of various computing tasks etc.}

\subsection{ \textbf{Abstract Continuum Computing Model  $AC^2M$ or $CM$ for short } } 

  The basic data objects and operations in this model are as follows :
\begin{itemize}

\item \textbf{Basic data objects}
\\ algebraic closure of meromorphic functions (over suitable domain)
\vspace{5pt}
\item \textbf{Basic unit operations}

\begin{itemize}
\item \vspace{5pt}\textbf{field operations}: $+, -, \times, \div$\\
and additionally for real quantities comparison $(<, =)$

\item \vspace{5pt}\textbf{analytic operations}
   \begin{itemize}
      \item \textbf{integration}: $\int_C f(z) dz$ this is a binary operation,
               function $f$ and specification contour $C$ are the two operands

        \item \textbf{differentiation}: $\frac{\partial f}{\partial z_i}$ is also a
               binary operation with inputs $f$ and $z_i$.
    \end{itemize}
\vspace{5pt}

\end{itemize}
\end{itemize}
(a word of caution when comparing with TM - there is no such thing as
``conservation of difficulty'' across the models.)

\newpage
\subsection{ \textbf{ Extension of $P \neq NP$ conjecture from TM theory}}
\vspace{5pt}
\begin{itemize}
\item Computing models can be organized in the same way as Cantor
  organized infinite sets i.e. according to the cardinal number
\\  of set of all possible data objects and machines in the model.
\vspace{5pt}
\item At each cardinal number in the sequence $\beth_0 < \beth_1 <
  \beth_2 < \ldots$, you have models of computation, and corresponding $P \neq NP$ question
\begin{eqnarray*}
  TM & \longrightarrow ~~~ \beth_0, & ~~~~  cardinal number (\mathbb{Q})
\\CM & \longrightarrow ~~~ \beth_1, & ~~~~  cardinal number (\mathbb{R})
\end{eqnarray*}
\vspace{0pt}
\item It appears that $P \neq NP$ problem for TM is just the first member in a
  sequence of strict inclusions
\[
P(TM) \subsetneq NP(TM) \subsetneq P(CM) \subsetneq NP(CM)
\subsetneq P(\beth_2) \ldots
\]
\vspace{0pt}
\item An interesting question across adjacent level:

\begin{tabular}{ll}
$NP(TM) \subsetneq P(CM)$
\end{tabular}
\vspace{5pt}
\\i.e., is non-determisitic computing at any one level no more powerful
than deterministic computing at the next level ?
\end{itemize}

\subsection{ \textbf{Our Research Program}}
Is aimed at understanding the following:
\vspace{5pt}
\begin{itemize}
\item how to construct physical machines supporting Continuum Computing ?
\vspace{5pt}
\item Can non-deterministic computing at the TM level be simulated by
deterministic computing at next level i.e. by Continuum Computing ? \\
\vspace{5pt}
 \textbf{and a harder question: } 
\item  to what extent can one approximate deterministic
  Continuum Computing by deterministic computing in TM?
\vspace{5pt}
\item initially, approximate cross-simulation was meant to be "stop-gap"
   measure, but now it appears that building continuum machine may take
   many years, hence simulation becomes more important.
\vspace{5pt}
 \item Floating Point numbers allows approximation of reals by rationals
\vspace{5pt}
 \item Similarly, one can approximate functions by other simpler functions
 \vspace{5pt}
 \item In this investigation, length of "binary encoding" of data object does not have the same fundamental significance as in TM theory. 
       Instead, other properties of the problem space seem more important
\end{itemize}

\newpage
\section{Extensions of proof systems for more efficient proofs of optimality and non-satisfiability }

\subsection{\textbf{Introduction} }
We are interested in converting results of continuum based algorithms so that proofs of optimality or non-existence of solutions can be verified on the standard model.

\vspace{5pt} Both involve proofs of non-negativity of functions.

\begin{itemize}

\item \textbf{Optimality}:

\vspace{5pt} Proving that $\textbf{x}_{min}$ is a global minimum of
$f(x_1, x_2, \ldots, x_n)$ is equivalent to showing that for $f_{min}
= f(\textbf{x}_{min})$, we have

\[
f(\textbf{x}) - f_{min} \ge 0 \hspace{15pt} \forall \textbf{x} \in \mathbb{R}^n
\]

\item \textbf{Non-satisfiability}:
Each variable $x_i = \pm 1 $ \\
let $V$ denote the variety defined by $x_i^2 = 1 \hspace{15pt}\forall i$ \\
For each clause associate a polynomial. e.g. \\
 $C = x_i \vee \bar{x_j} \vee x_k$, associated
$p_c(\textbf{x}) = \left[1 - x_i\right] \cdot \left[1 + x_j \right] \cdot
\left[1 - x_k\right]$\\

For a $\pm 1$ vector \textbf{x}, \hspace{5pt}
%\[
$p_c (\textbf{x})=
\begin{cases}
0 & \text{if}~\textbf{x}~\text{is a satisfying assignment}\\
8 & \text{otherwise}
\end{cases}$\\
%\]
Define $f = \sum_{C \in \text{Clauses}} p_c(\textbf{x})$.

\vspace{3pt}Then $f(\textbf{x}) > 0 \hspace{15pt} \forall \textbf{x} \in V $
gives proof of non-satisfiability.

\end{itemize}

\subsection{\textbf{Relation of Positivity Proofs to Hilbert's $17^{\text{th}}$ problem}}
\subsubsection{\textbf{introduction}}
\begin{itemize}
\item One approach to proving positivity of a polynomial is
  to express it as sums of squares of other functions

\item \textbf{Using polynomials}:

\textbf{Hilbert} (1888) realized that it is \textbf{not} always possible ,
but gave a \textbf{non-constructive} proof of the existence of a counter
example.

\item \textbf{Using rational functions}:

Is always possible -- \textbf{Artin's} (1926)  solution of Hilbert's 17th problem.\\
  But exponentially many terms of high degree required (\textbf{Pfister}).

\item \textbf{Concrete examples} of Hilbert's result took long time to
  construct.

\hspace{5pt} \textbf{First counter example - Motzkin}
\[
f(x, y, z) = z^6 + x^4 y^2 + x^2 y^4 - 3 x^2 y^2 z^2
\]

\hspace{5pt} \textbf{Further examples - Lam, Robinson} and others.

\item  \textbf{Computational results}:

 For these ``counterexamples'', we have computed expressions as
sums of squares of polynomials, with a modified interpretation. \\
Algorithm underlying the solver was described in previous lecture
[Karmarkar, MIT91, IPCO92] 
\end{itemize}
\newpage
\subsubsection{\textbf{Computational approach to positivity proofs}}
\begin{itemize}
\item \vspace{-3pt} Consider a homogeneous polynomial $f(x_1, x_2,
  \ldots, x_n)$ of degree $d$ in $n$ real variables $x_1, x_2, \ldots,
  x_n$, which is non-negative everywhere.

\item \vspace{5pt} To show that a function is non-negative on a
  compact set, it is enough to show that the function is
  \textbf{non-negative} at all its \textbf{critical points}.

(note we have homogeneous polynomials over projective space).
Therefore non negativity at set of infinite number of points in \bbbr $^n$  is implied by non-negativity at a finite subset of points viz the critical points.
\item  \vspace{5pt} This can be achieved if we \\

 \indent \indent (1) construct a \textbf{variety containing
  all critical points} of the function and\\
\vspace{3pt} \indent \indent (2) construct an expression of the
function as sums of squares

\hspace{11pt} of polynomials \textbf{in the coordinate
  ring of that variety}, instead of

\hspace{11pt} the polynomial ring $\Re(x_1, \ldots, x_n)$.

\item \vspace{5pt} Additionally, without loss of generality, we impose a spherical
  constraint
%\[
\begin{equation}
x_1^2 + x_2^2 + \ldots + x_n^2 = 1
\end{equation}
%\]

\item Each point in the real projective space is covered twice 
 (by a pair of antipodal points) in this representation.

\item \vspace{5pt}  A critical point of the function satisfies \\

\begin{equation}
\frac{\partial f}{\partial x_i} = \lambda x_i, \ \ i = 1, ... ,n
\end{equation}

where $\lambda$ is the Lagrange multiplier.

\end{itemize}
%\section{Computational approach to positivity proofs (Contd.)}
\begin{itemize}

\item \vspace{5pt} We now work in $\Re^{n+1}$, and points in this
  expanded space will be denoted by $(x_1,x_2,...,x_n,\lambda)$

\item \vspace{5pt} The equations (1) and (2) define an algebraic
  variety $U$ (i.e an algebraic set -- in our terminology a variety does
  not have to be irreducible).

\item \vspace{10pt} Our approach is to construct a variety $V$ such
  that $U\subseteq V \subseteq \Re^{n+1} $\\ and in the co-ordinate
  ring of $V$, construct an expression for f as sums of squares

\begin{center}

$\sum_i S_i ^2 (x_1, x_2, \ldots, x_n, \lambda)$

\end{center}

\item \vspace{5pt} We also produce an explicit expression for \\

\begin{center}
\begin{tabular}{ll}

\vspace{3pt}& $ f\ -\ \sum_i S_i^2\ (x_1, x_2, \ldots, x_n, \lambda)\ $ \\
\vspace{5pt} as &
$\sum_j\ a_j\ (x,\lambda)\ g_j\ (x_1, x_2, \ldots, x_n, \lambda)\ $ \\
\vspace{3pt} where & $g_j\ (x_1, x_2, \ldots, x_n, \lambda)\ = \ 0 $ \\

\end{tabular}
\end{center}

 are the defining equations for $V$. \\

\end{itemize}
\subsubsection{\textbf{Samples of computer generated proofs}}
For small examples, the proofs are short enough to be displayed on viewgraphs : \\
\begin{footnotesize}\\
\textbf{Motzkin's first counter example} \\
\begin{eqnarray*}
%\begin{array}{lcr}
f(x,y,z) & = & z^6 + x^4y^2 + x^2y^4 - 3x^2y^2z^2 \nonumber \\
f & \equiv & (\frac{1}{4}S_1^2 + S_2^2 + S_3^2 + S_4^2 + \frac{3}{4}S_5^2)\ \text{mod}~ V \nonumber \\
S_1 & = & xy(x^2-y^2)\nonumber \\
S_2 & = & x^4 + y^4 - 2x^2 - 2y^2 + x^2y^2 + 1 \nonumber \\
S_3 & = & xz (x^2 + 2y^2 - 1) \nonumber \\
S_4 & = & yz (2x^2 + y^2 - 1) \nonumber\\
S_5 & = & xy(3x^2 + 3y^2 - 2) \nonumber\\
%\end{array}
\end{eqnarray*}

\vspace{-3pt}\textbf{Robinson's counter example} \\
\begin{eqnarray*}
f & = & x^6 + y^6 + z^6 - (x^4y^2 + x^2y^4 + y^4z^2 + y^2z^4 + z^4x^2 + z^2x^4 ) + 3x^2y^2z^2     \nonumber \\
f & \equiv & (S_1^2 + \frac{3}{4}S_2^2 + \frac{1}{4}S_3^2 + S_4^2 + S_5^2)\ \text{mod}~ V \nonumber \\
S_1 & = & -x^3y + xy^3 \nonumber \\
S_2 & = & -1 + 3x^2 - 2x^4 -4x^2y^2 + 2y^2 \nonumber \\
S_3 & = & 1 - x^2 - 2x^4 + 4x^2y^2 -4y^2 + 4y^4 \nonumber \\
S_4 & = & -2x^3z - xy^2z + xz \nonumber \\
S_5 & = & -x^2yz - 2y^3z + yz \nonumber \\
\end{eqnarray*}
\end{footnotesize}

%\vspace{80pt}
%\begin{center}
%\vspace{10pt}
%$\Large{Extensions of previous proof systems}
%\vspace{10pt}
%\end{center}

\subsection{Rules for positivity proofs}
\subsubsection{\textbf{Basic Rules}}
%\setbeamertemplate{enumerate items}[default]
\begin{enumerate}

\item  \textbf{Starting primitives}:
\begin{enumerate}
\item \vspace{5pt} \textbf{Constant functions}

\vspace{3pt} Let $\alpha \in \mathbb{R}_+$ be a positive scalar.

\vspace{3pt} If $f(x) = \alpha$, then $f > 0$.

\item \vspace{5pt} \textbf{Square functions}

\vspace{3pt} If $f(x) = g^2(x)$, then $f \ge 0$.
\end{enumerate}
\item \textbf{Algebraic operations preserving positivity}:

\begin{center}
\vspace{2pt} ~\,If $f \ge 0$ and $ g \ge 0$, then $f + g \ge 0$

\vspace{5pt} If $f \ge 0$ and $ g \ge 0$, then $f \cdot g \ge 0$
\end{center}
\item  \textbf{Positivity restricted to variety defining constraints}:

Suppose there are integrality constraints, $x_i^2 = 1$, or other
constraints like $A\textbf{x} = \textbf{b}$ etc. They define an algebraic
set or subvariety $V$ of $\mathbb{R}^n$. It is understood that algebraic operations are
to be interpreted upto equivalence classes of functions on $V$.
\end{enumerate}

\vspace{3pt} There has been extensive exploration of application of rules of this type along with rules for cuts based on integrality constraints. see e.g. Lova\'sz, Schrijver, Grigoriev,Chlamatac,Schoenebeck, Tulsiani,Worah, and references therein

\subsubsection{\textbf{New rules}}
We strengthen the previous proof systems by adding the following new rules.
\begin{enumerate}
\item \textbf{Substitution or Composition rule}:

\vspace{3pt} Let $f(x_1, \ldots, x_n)$ and $g_i(y_1, \ldots, y_m)$ for
$1 \leq i \leq n$ be real polynomials 

\vspace{3pt} Let $h(y_1, \ldots, y_m) = f(g_1(y_1, \ldots, y_{m}),
\ldots, g_n(y_1, \ldots, y_m))$.
\begin{enumerate}
\item \vspace{5pt} If $f(x_1, \ldots, x_n) \ge 0$ for all $\textbf{x} \in
  \mathbb{R}^n$, then $h(y_1, \ldots, y_m) \ge 0$.
\item \vspace{5pt} If $f(x_1, \ldots, x_n) \ge 0$ for all
  $\textbf{x} \in \mathbb{R}^n_+$ and $g_i(y_1, \ldots, y_m) \ge 0$ for
  all $\textbf{y} \in \mathbb{R}^m$, then $h(y_1, \ldots, y_m) \ge 0$.

\end{enumerate}

\vspace{5pt}
\item \textbf{Division rule}:
In some sense, this rule is not new, since Artin's approach uses rational functions.
However, number of rational terms in that approach can be exponential, but having a separate rule enables selective application for more efficient proof. \\
Suppose $g(x) > 0, f(x) \ge 0$ and $g(x) \mid f(x)$.
\begin{center}
If $f(x) = g(x) \cdot h(x)$, then $h(x) \ge 0$.
\end{center}
This rule reduces the degree.

\vspace{5pt}
\item \textbf{Odd Radical rule}:

Suppose $f(\textbf{x})$ is a perfect odd ($(2k+1)^{\text{th}}$)
power. Let $g(\textbf{x})$ be the $(2k+1)^{th}$ root of $f(\textbf{x})$,
i.e.  $f(\textbf{x}) = g^{2k+1}(\textbf{x})$.
\begin{center}
If $f(\textbf{x}) \ge 0$, then $g(\textbf{x}) \ge 0$.
\end{center}
This is one of the rare rules that reduces the degree.\\

 Notes:
 \begin{enumerate}
   \item Certain inequalities ( e.g. Cauchy-Schwartz ) are proved symbolically, and have proof-length O(1), independent of n. When invoking such inequalities during application of the composition rule, only length involved in specifying the substitution should be counted as addition to the length of the proof.      
   \item Lower bounds based on proof system without these rules don't show intrinsic difficulty of any problem, but only show intrinsic limitations of the proof system used
 \end{enumerate}
  
\end{enumerate}

\section{\textbf{Maximum Independent Set}}
\subsection{\textbf{Introduction}}
 \begin{itemize}
 \item Let $G = (V, E)$ be a graph with vertex set $V$ and edge set $E$.
 \item A subset $U \subseteq V$ of nodes is called \emph{independent}
%   if for every pair of nodes $i, j \in U$, we have $(i, j) \notin E$.
    if for every edge $(i, j) \in E$, either  $i \notin U$  or $j \notin U$
 \item Problem: Find the largest independent set.
 \item $\pm 1$ integer programming formulation

 Constraints:

 \begin{center}
 \begin{tabular}{rlrl}
 $w_i $ & = & $\pm 1;$ & $\,i = 1, \ldots, n$\\
 $w_i + w_j$ & $\leq$ & 0; & $\,\forall (i, j) \in E$\\
 \end{tabular}
 \end{center}
 Objective: Maximize $\sum w_i$
 \item Relaxation of the feasible set allows $-1 \leq w_i \leq 1$
 \item Each extreme point of an LP has co-ordinates
   that are $\pm 1$ or 0.
 \item  For a $\pm 1$ solution, $\sum_{i} w_i^2 = n$.
 \item Non-convex objective: $\sum w_i + \beta \sum w_i^2$
 \end{itemize}
 \subsection{\textbf{Defining a sharper polytope}}

 \vspace{10pt}
 There are additional inequalities
 \begin{itemize}
 \item \vspace{5pt} that hold for the solution points, i.e. feasible
   points with $\pm 1$ solution
 \item \vspace{5pt} but cannot be expressed as non-negative linear
   combinations of basic inequalities
 \item \vspace{5pt} number of such inequalities grows exponentially with $n$
 \item \vspace{5pt} significant effort in polyhedral combinatorics is
   aimed at identifying such constraints (e.g. see Schrijver\cite{Schrijver}).
 \end{itemize}

 \subsection{\textbf{Example: odd cycle inequalities}}
 \begin{itemize}
 \item Consider an odd cycle in the graph, having $2k + 1$ vertices. At
   most $k$ of the $2k + 1$ vertices can be in an independent set.
 \begin{figure}
 \includegraphics[scale=0.4]{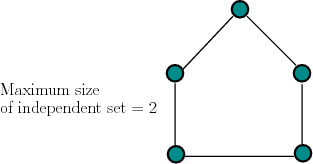}
 \end{figure}
% %%%% A figure is due here.
 \item Let $i_1, i_2, \ldots, i_{2k+1}$ be the vertices in the cycle;
   hence $w_{i_1} + w_{i_2} + \ldots + w_{i_{2k+1}} \leq -1$.
 \item This inequality is sharper, since the previous inequalities only imply :\\
   $w_{i_1} + w_{i_2} + \ldots + w_{i_{2k+1}} \leq 0$.
 \item Number of odd cycles, and hence the corresponding number of
   inequalities, can grow exponentially with $n$.
 \end{itemize}
 \subsection{\textbf{Inequalities for other subgraphs}}
 \begin{itemize}
 \item Similar to odd cycles, there are constraints based on other
   types of subgraphs.
 \item Let $$\sum_{i \in G_0} a_i \cdot w_i \leq b$$ be a constraint
   based on a fixed graph $G_0$.
 \item Let $\widetilde{G}$ be a graph obtained by an odd sub-division
   of the edges of $G_0$, i.e. replace an edge by an odd path.
 \vspace{5pt}
 \begin{figure}
 \includegraphics[scale=0.4]{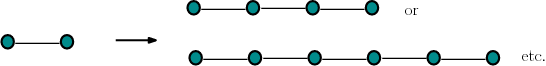}
 \end{figure}
%% %%%% figure is due here....
 \item \vspace{5pt} Then sum of newly added vertices is always
   non-positive, i.e. $\sum_{i \in \text{New vertex}} w_i \leq 0$.
 \end{itemize}
 \vspace{-5pt}
 \begin{itemize}
 \item Furthermore, to achieve the equality,
     $\sum_{i \in \text{New vertex}} w_i = 0$,there are \\
    only two ways to assign values to new vertices
   and each corresponds to a particular feasible assignment to
   the end points of the original edge.
 \vspace{5pt}
 \begin{figure}
 \includegraphics[scale=0.25]{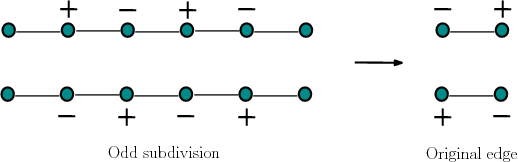}
 \end{figure}
 \item Using this observation, it is easy to show that the inequality
   for $G_0$ implies a similar inequality for $\widetilde{G}$, i.e.

 \[
 \sum_{i \in G_0} a_i \cdot w_i \leq b ~~~~\Rightarrow~~~~  \sum_{i \in \widetilde{G} }
 \widetilde{a_i} \cdot w_i \leq b
 \]
 \item In this context, we are using a more restricted notion of
   homeomorphism.
 \item Specifically, in order for two paths to be homeomorphic \\
   both must have odd length or both must have even length.
 \item We call this ``parity-respecting'' homeomorphisms, and denote it by
 \vspace{-5pt}
 \[
 \widetilde{G} \cong G_0 (\text{mod} 2)
 \]
  (different authors use different
   terminology/notation, e.g. see [])
 \end{itemize}
 \begin{itemize}
 \item Given a fixed graph $G_0$ and an input graphs $G$, \\
   an inequality for $G_0 \Rightarrow $ an inequality for each subgraph $\widetilde{G}$ of $G$
   that is homeomorphic to $G_0$.
 \item \vspace{3pt} Number of such subgraphs can grow exponentially
   with $n$.
 \item \vspace{3pt} As an example, let $G_0$ be a triangle and $G$ be a given graph. \\
   Then the subgraphs of $G$ that are homeomorphic to $G_0$ respecting parity, \\
   are exactly the odd cycles of $G$. \\
    All odd cycle inequalities for $G$ are obtained from inequalities for $G_0$.
 \item \vspace{3pt} We will use this example to show how the combined
   effect of all odd cycle inequalities can be computed in polynomial
   time in the continuum model.
 \item \vspace{3pt} Observe that all the inequalities in this example
   are linear, which simplifies our exposition.
 \item \vspace{3pt} However, these techniques are not limited to the
   linear case. Later, we will show an example of non-linear,
   non-convex inequalities as well.
 \end{itemize}
 
 \newpage
 \section{\textbf{Combining Effect of Exponential Number of Inequalities}}
 
 \subsection{\textbf{Projectively Invariant Metric}}
 \begin{itemize}
 \item Consider the projectively invariant metric we use in linear
   programming algorithm. Let
 \[
 \triangle = \left\lbrace x_i ~\vert~ x_i \geq 0, \sum x_i =1
 \right\rbrace
 \]

 be the simplex used in the algorithm.
 \item Let $\textbf{x}, \textbf{y} \ \epsilon \ int\left( \triangle \right) $
   be two interior points in the simplex.  \\
   Projectively invariant distance between $\textbf{x}$ and $\textbf{y}$, based on $p$-norms :
 \begin{align*}
 d(\textbf{x}, \textbf{y}) = \frac{1}{2} \cdot \Bigg(\sum_i
 \left(\frac{x_i}{y_i}-\frac{y_i}{x_i}\right)^p \Bigg)^{\frac{1}{p}}
 \end{align*}
 \item The infinitesimal version of d gives the Riemannian metric
   $g(\textbf{x})$ for $p=2$ and Riemann-Finsler metric for $p > 2$.
 \begin{align*}
 g^p _{ij} (x)dx_idx_j = \sum_i \left( \dfrac{dx_i}{x_i}\right) ^p
 \end{align*}

 The performance of the algorithm depends on curvature in this metric.
Karmarkar \cite{Kar4}
 
 \end{itemize}
 
 \begin{itemize}
 \item For inequalities in the form $A \textbf{x} \leq b$, let $\textbf{x},
   \textbf{y}$ be two interior points, and let $\textbf{s}, \textbf{t}$ be the
   corresponding slack variables i.e.
 \[
 A\textbf{x} + \textbf{s} = \textbf{b}, ~~~~~ A\textbf{y} + \textbf{t} = \textbf{b}
 \]

 Then projectively invariant distance between $\textbf{x}$ and $\textbf{y}$
 is given by
 \begin{align*}
 d^p(\textbf{x}, \textbf{y}) = \sum_i \left[\frac{1}{2} \cdot
   \left(\frac{s_i}{t_i}-\frac{t_i}{s_i}\right)\right]^p
 % d(\textbf{x}, \textbf{y}) = \sum_i \|\dfrac{1}{2} \left[
 %   \dfrac{s_i}{t_i}-\dfrac{t_i}{s_i}\right] \| _p
 \end{align*}
 \item Let $\textbf{s}_0 = $ slack variable for the current interior point
   $\textbf{x}_0$ (constant) \\and $\textbf{s} = $ slack variable for the next
   (unknown) interior point $\textbf{x}$.%% , to be found as the output of the
   %% current step.

 \item Observe the particular \textbf{``distributive``} form of the function
   $d^p (\textbf{x}, \textbf{x}_0)$.

 If $S$ is the set of all slack variables, then $d^p$ distributes
 linearly over $S$.
 \begin{align*}
 d^p(\textbf{x}, \textbf{x}_0) = \sum_{s \epsilon S} \psi(s)
 \end{align*}

 where the $\psi(s)$ for the \textbf{individual slack variable} is
 \begin{eqnarray*}
 \psi(s) & =  & \left[\frac{1}{2} \cdot
   \left(\frac{s}{s_0}-\frac{s_0}{s}\right)\right]^p
 \end{eqnarray*}

 \end{itemize}
 \subsection{\textbf{Computing exponential sums efficiently}}
 \begin{itemize}
 \item Note that the potential function $\psi$ used in the interior
   point methods also has the same \textbf{distributive} form, \\
   where $\psi$
   could be a non-linear rational or transcendental function of
   individual slack variable.
   \vspace{2pt}
 \item We are interested in the case when $S$ is exponentially large.
   \vspace{2pt}
 \item Such sums can be evaluated efficiently for \\
   \begin{itemize}
   \item Large class of problems involving exponential number of inequalities, \\
   \vspace{4pt}
   \item when $\psi(s)$ belongs to certain special parametric family of functions
   such as exponential function, e.g $\psi (s) = e ^{-zs}$.
   \end{itemize}
 \vspace{2pt}
 \item While the actual function $\psi(s)$ of interest is not of this form,\\
    it can be expressed as a linear superposition of functions of this form.
 \vspace{2pt}
 \item E.g. techniques such as \\
   Laplace transform or Fourier transform and their inverse transforms  enable\\
   expressing $\psi$ as (infinite)superposition of (real or complex) exponentials.
 \vspace{2pt}
 \item For approximate cross-simulation on the standard model, \\
   we use a suitable \textbf{finite superposition.}
 \end{itemize}
 
 \section{\textbf{Combined Effect of Inequalities for all Odd Cycles in a
   Graph}} 
\subsection{\textbf{introduction}}
   Contribution of an individual slack variable $s$: $\psi (s) =
 e ^{-zs}$ where $z \epsilon ~\mathbb{C}$.\\
\vspace{4pt}
 \textbf{Goal}: We want to find closed form meromorphic function for \\
 total contribution of all slack variables.\\
\vspace{4pt}
 \textbf{Edge Matrix}
 \begin{itemize}
 \item For a single edge, we have $\dfrac{w_i + w_j}{2} \leq 0$ \\
   Note: since each node in a cycle has degree 2, we are dividing by 2\\
   For more general graph minors, the weighting factors are different.

 \item Introducing slack variable $s_{ij},$ we get
 $\dfrac{w_i + w_j}{2} + s_{ij} =0 $. Then,
 \begin{eqnarray*}
 s_{ij} &  = & -\dfrac{w_i + w_j}{2}\\
  \psi & = & e^{-z \cdot s_{ij}}=e^{z\left\lbrace \frac{w_i + w_j}{2} \right\rbrace }
 \end{eqnarray*}
 \item Define edge matrix $A(z, \textbf{w})$ over the field of meromorphic
   functions as
 %\begin{align*}
 \[
   A_{ij}(z, \textbf{w}) =
 \begin{cases}
     e^{\frac{z}{2}\left\lbrace w_i + w_j \right\rbrace } & \text{if $(i, j)$ is an edge}\\
     0 & \text{otherwise.}
   \end{cases}
 \]

 \end{itemize}
 \begin{itemize}
 \item For a single slack variable $s$, there is a relation between the
   derivative operators $\dfrac{\partial}{\partial s}$ and
   $\dfrac{\partial}{\partial z}$
 \begin{align*}
 s \dfrac{\partial \psi}{\partial s} = z \dfrac{\partial \psi}{\partial z}
 \end{align*}
 \item For the edge matrix $A(z, \textbf{w})$, derivative with respect to
   $w_i,$ the co-ordinates of the interior point, are expressed in
   terms of the ``J-products''.
 \begin{align*}
 \dfrac{\partial A}{\partial w_p} = \frac{z}{2} \left\lbrace I
 \bigcirc \hspace{-0.9em} {\raisebox{0.00ex} {\footnotesize J}} \hspace{-0.3em} \raisebox{-1.50ex} {\tiny p}
 \hspace{0.7em}A + A
 \bigcirc \hspace{-0.9em} {\raisebox{0.00ex} {\footnotesize J}} \hspace{-0.3em} \raisebox{-1.50ex} {\tiny p}
 \hspace{0.7em}I \right\rbrace
 \end{align*}
 \item This gives an autonomous differential equation for $A$, \\
   which allows us to create recurrence relation connecting \\
  \begin{itemize}
   \item higher powers of $A$ to lower powers and
   \item \vspace{3pt} higher derivatives of $A$ to lower order derivatives.
 \end{itemize}
 \item  These recurrence relations lead to closed form expressions for higher derivatives and higher powers.
 \end{itemize}
 
 \subsection{\textbf{Permitting additional inequalities}}
 %\begin{itemize}
 %\item

 We are interested in all odd cycles in the graph.

 \vspace{5pt} A odd cycle is a special case of a closed odd walk.\\

 \vspace{5pt} Since a closed odd walk contains an odd cycle as a
 subgraph, a sharper version of inequality is also valid for closed odd
 walk.

 \vspace{5pt} However this inequality is implied by the odd cycle
 inequalities.

 \vspace{5pt} Therefore, from the point of view of formulation, these
 additional inequalities are superfluous:

 \vspace{5pt} useless $\rightarrow$ since they are implied by other
 inequalities but \\ harmless $\rightarrow$ they are still valid.

 \vspace{5pt} However from the point of view of obtaining a closed form
 expression that can be evaluated in polynomial time, they are
 essential.

 \section{\textbf{Getting the closed form expression for the combined effect
   of inequalities}}
\subsection{\textbf{Introduction}}
 \textbf{Steps}:

 \vspace{10pt}(1) Get expressions for the effect of the \textbf{implied}
 inequalities for the following:

 \vspace{5pt}~~ (a) For a \textbf{given walk} $W$ of \textbf{length $l$} ~:~ $\psi_W$

 \vspace{5pt}~~ (b) For a \textbf{given pair} of nodes $i$ and $j$ , \textbf{all} walks of length $l$ ~:~ $\psi_{i j}$

 \vspace{5pt}~~ (c) In the \textbf{entire} graph, all \textbf{closed} walks of length $l$  ~:~ $\psi_l$

 \vspace{5pt}~~ (d) In the entire graph, all closed \textbf{odd} walks of \textbf{length at most $n$}  ~:~ $\psi$

 \vspace{10pt} (2) Transitioning from expressions for implied
 inequalities to expressions for \textbf{sharper} inequalities for the
 following:

 \vspace{5pt}~~ (a) All closed odd walks of length at most $n$ in the
 entire graph ~:~ $\widetilde{\psi}$

 %\item $A \underset{p}{\JP} B $, $C  \underset{2}{\otimes} D$

 \subsection{\textbf{Closed form expression: Step (1a)}\\}

 \textbf{Single walk $W$ of length $l$}:

 \vspace{5pt}Consider a walk $W$ with $l$ edges $i_1, i_2, \ldots,
 i_{l+1}$. Summing the edge inequalities, we have

 \[
 \frac{1}{2}w_1 + w_2 + \ldots + w_l + \frac{1}{2}w_{l+1} \leq 0
 \]

 Let $s_w$ be a slack variable for the ``implied'' inequality for walk
 $W$.
 \begin{eqnarray*}
 s_w & = & \sum_{e \in W} s_e\nonumber\\
 \psi_W(s_w) & = & e^{-z \cdot s_w} = e^{-z \sum_{e \in W} s_e} = \prod_{e \in W} e^{-z \cdot s_e}\nonumber\\
 \psi_W & = & A_{i_1 i_2} A_{i_2 i_3} \cdots A_{i_l i_{l+1}}
 \end{eqnarray*}

 %\end{itemize}
 \subsection{\textbf{Closed form expression: Steps (1b) \& (1c)}\\}
  \textbf{Between given pair of nodes $i$ and $j$ effect of all walks of length} $l$   :

 \vspace{5pt} Let $i_1 = i$ and $i_{l+1} = j$.
 \begin{eqnarray*}
 \psi_{i j} & = & \sum_{i_2, i_3, \ldots, i_l} A_{i i_2} A_{i_2 i_3} \cdots A_{i_l j}\nonumber\\
 \psi_{i j} & = & A^l_{i j}\nonumber
 \end{eqnarray*}

 \textbf{All closed walks of length $l$ in the entire graph}:

 \[
 \psi_l = \sum_{i} A_{i i}^{l} = tr \{A^{l}(z, \textbf{w})\}
 \]
 However, a closed walk of length $l$ is counted $2l$ times in the
 expression above,\\
 \begin{itemize}
\item Due to $l$ vertices of the walk and \\
\item The two senses of traversal along the walk.\\
\end{itemize}

 Compensating for this repetition, we have

 \[
 \psi_l = \frac{tr\{A^{l}(z, \textbf{w})\}}{2l}
 \]

 \subsection{\textbf{Closed form expression: Step (1d)}\\}

 \textbf{All closed odd walks of length at most $n$ in the entire graph}:

 \vspace{3pt} Let $k_{max} = \lfloor \frac{n-1}{2} \rfloor$.
 \begin{eqnarray*}
 \psi & = & \sum_{k = 1}^{k_{max}} \psi_{2k+1}\nonumber\\
 & = & \sum_{k = 1}^{k_{max}} \frac{1}{2 \cdot (2k+1) } tr\{A^{2k+1}(z, \textbf{w})\}\nonumber\\
 & = & tr\{B(z,\mathbf{ w})\}\nonumber
 \end{eqnarray*}
 where
 \[
 B(z,\mathbf{w}) = \sum_{k = 1}^{k_{max}} \frac{A^{2k+1}}{2 \cdot (2k+1)}
 \]

 \subsection{\textbf{Transitioning from effect of implied inequalities to
   sharper ones}\\}

 Let $i_1, i_2, \ldots, i_{l+1}$ be a walk of length $l$; its implied
 inequality is given by
 \[
 \frac{w_{i_1}}{2} + w_{i_2} + \ldots + w_{i_l} + \frac{w_{i_{l+1}}}{2} \leq 0
 \]

 If $i_{l+1} = i_1$, we have a closed walk with the implied inequality
 given by
 \[
 w_{i_1} + w_{i_2} + \ldots + w_{i_l} \leq 0
 \]

 For closed odd walks ( $l = 2k+1 $), we get the following
 \textbf{sharper} inequality.
 \[
 w_{i_1} + w_{i_2} + \ldots + w_{i_l} \leq -1
 \]

 since at most $k$ nodes can be in an independent set.

 If $s$ is the slack variable for the implied inequality, let
 $\widetilde{s}$ denote the slack variable for the corresponding
 sharper inequality.
 \begin{eqnarray*}
 \widetilde{s} & = & s -1 \nonumber\\
 e^{-z \cdot \widetilde{s}} & = & e^{z} \cdot e^{-s}\nonumber\\
 \widetilde{\psi_l} & = & e^{z} \psi_l\nonumber\\
 \end{eqnarray*}

 \subsection{\textbf{Closed form expression: Step (2a)}\\}

 Transitioning from the expression for implied inequality for a closed
 odd walk to the sharper inequality  involves multiplication by
 $e^z$, which is independent of the length of the walk.

 Therefore, for $l$ odd, $\widetilde{\psi_l}$ is given by
 \[
 \widetilde{\psi_l} = e^z \cdot \frac{tr\{A^{l}(z, \textbf{w})\}}{2l}
 \]

 And $\widetilde{\psi}$ is given by
 \begin{eqnarray*}
 \widetilde{\psi} & = & e^z \cdot B(z, \textbf{w})%% \psi\nonumber\\
 %% & = & e^z \cdot tr\left\{\sum_{k = 1}^{k_{max}} \frac{A^{2k+1}}{2k+1}\right\}
 \end{eqnarray*}

 %% Since $e^z$ is a multiplicative factor that is independent of the
 %% length of odd cycles, the closed form expression for
 %% $\widetilde{\psi}$ is as follows:

 %% \[
 %% \widetilde{\psi} = e^z B(z, \textbf{w}) (= \widetilde{B}(z, \textbf{w}))
 %% \]

 It is trivial to see that this can be evaluated in polynomial time in
 the continuum model.
 \begin{itemize}
\item   Straightforward  evaluation of $\widetilde{\psi}$ as written, will take $O(n^4)$ operations. But with some rearrangement,  it is possible to evaluate it with $O(n^3)$ operations. 
\paragraph{    We also need first two derivatives of $\psi$ for constructing the trajectory in the main algorithm. We need the concept of J-products to evaluate derivatives efficiently.}

  \end{itemize}
\newpage
 \section{\textbf{Join Product}}
 \subsection{\textbf{Introduction}}
We split the tensor product involving contraction into two steps
The first step is similar to join operation on relations in a relational data base\\
\begin{figure}
 \includegraphics[scale=0.4]{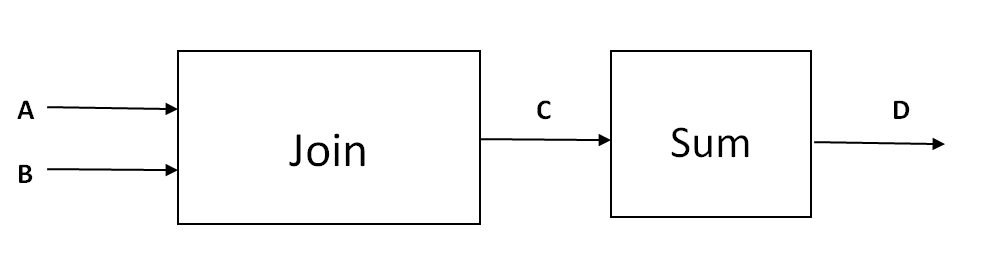}
\end{figure}
\textbf{Simple examples of J-products --}\\
\begin{enumerate}
\item \textbf{ Matrix Multiplication :}
\begin{align*}
C = A.B, \hspace{5pt} C_{pr} = \sum_{q} A_{pq} B_{qr}
\end{align*}
corresponding "join" product (definition):
\begin{align*}
C_{pqr} = A_{pq} B_{qr}
\end{align*}
without the summation (contraction) over q \\
corresponding "join" product (notation):
\begin{align*}
C = A
\bigcirc \hspace{-0.9em} {\raisebox{0.00ex} {\footnotesize J}} \hspace{-0.7em} \raisebox{-1.70ex}{ \tiny q}
\hspace{0.5em} B
\end{align*}

\end{enumerate}

\begin{enumerate}

\item \textbf{Dot Product of vectors :}
\begin{align*}
\textbf{c = a.b} = \sum_i a_i b^i
\end{align*}

Corresponding join product (definition) :
\begin{align*}
   c_i = a_ib^i
\end{align*}

Corresponding join product (notation) :
\begin{align*}
\textbf{c} = \textbf{a} \hspace{0.5em}  \hspace{-0.5em}
\bigcirc \hspace{-0.9em} {\raisebox{0.00ex} {\footnotesize J}} \hspace{-0.5em} \raisebox{-1.70ex}{ \tiny i}
\hspace{0.9em} \textbf{b}
\end{align*}
\item \textbf{Hadamard Product of matrices :}\\
 is an example of J-product with two repeated indices\\
\begin{align*}
A_{pq} = B_{pq}C_{pq}
\end{align*}
\begin{align*}
A = B \hspace{0.5em}  \hspace{-0.5em}
\bigcirc \hspace{-0.9em} {\raisebox{0.00ex} {\footnotesize J}} \hspace{-1.0em} \raisebox{-1.50ex}{ \tiny p, \tiny q}
\hspace{0.7em} C
\end{align*}
\end{enumerate}
\subsection{\textbf{Tensor contraction}\\}
   \textbf{input tensors :}   A, B, with repeated covariant and contravarient indices , and \\
 \textbf{output tensor :}  C with summation over repeated indices $r_1 , r_2, ... r_m$.\\
\vspace{0pt}
\begin{align*}
C_{q_1q_2...q_lt_1t_2...t_u}^{p_1p_2...p_ks_1s_2...s_n} = \sum_{r_1r_2...r_m}A_{q_1q_2...q_lr_1r_2...r_m}^{p_1p_2...p_k}B_{t_1t_2...t_u}^{r_1r_2...r_ms_1s_2...s_n}
\end{align*}
\vspace{0pt}
If there are m repeated indices, rank(C)   =   rank(A)   +   rank(B)   -   2.m\\
\vspace{0pt}
Corresponding \textbf{join operator :}\\
\begin{itemize}
\item There is no summation over repeated indices \\
\item One copy of such indices is present in the output, enclosed in round () \\
\item They don't have significance as tensor indices but simply denote \\
an indexed family of tensors of the same rank as above \\
\end{itemize}
\vspace{0pt}
\begin{align*}
             [C_{q_1q_2...q_lt_1t_2...t_u}^{p_1p_2...p_ks_1s_2...s_n}]_({r_1r_2...r_m}) = A_{q_1q_2...q_lr_1r_2...r_m}^{p_1p_2...p_k}\hspace{0.5em}  \hspace{-0.5em} \bigcirc \hspace{-0.9em} {\raisebox{0.00ex} {\footnotesize J}}  \hspace{0.7em}B_{t_1t_2...t_u}^{r_1r_2...r_ms_1s_2...s_n}
\end{align*}

\vspace{0pt}

\textbf{Notation :}
\begin{align*}
C = A\hspace{-0.0em} \bigcirc \hspace{-0.9em} {\raisebox{0.00ex} {\footnotesize J}}  \hspace{0.7em}B
\end{align*}
Sometimes, the repeated indices are noted below $\hspace{0.5em} \bigcirc \hspace{-0.7em} {\raisebox{0.00ex} {\footnotesize J}} \hspace{0.5em} $ as shown :\\

\begin{align*}
C = A\hspace{1.0em} \bigcirc \hspace{-0.9em} {\raisebox{0.00ex} {\footnotesize J}} \hspace{-1.5em} \raisebox{-1.50ex}{ \tiny r1r2...rm} \hspace{0.7em}B
\end{align*}

\textbf{ Tensor Product as composition of join and summation operation :} \\
\begin{align*}
A.B = \sum_{r_1r_2...r_m}A\hspace{1.0em} \bigcirc \hspace{-0.9em} {\raisebox{0.00ex} {\footnotesize J}} \hspace{-1.5em} \raisebox{-1.50ex}{ \tiny r1r2...rm} \hspace{0.7em}B
\end{align*}

\subsection{\textbf{Properties of join product}}
\begin{itemize}
\item \textbf{Linear in both arguments :}\\
If $\alpha$   and $\beta$  are scalars\\
\begin{align*}
\{\alpha A + \beta B\} \bigcirc \hspace{-0.9em} {\raisebox{0.00ex} {\footnotesize J}}  \hspace{0.7em}C =
\alpha (A  \bigcirc \hspace{-0.9em} {\raisebox{0.00ex} {\footnotesize J}}  \hspace{0.7em} C) +
\beta (B  \bigcirc \hspace{-0.9em} {\raisebox{0.00ex} {\footnotesize J}}  \hspace{0.7em} C)
\end{align*}
If A, B, X are matrices of compatible dimensions,\\
\begin{align*}
X (A  \bigcirc \hspace{-0.9em} {\raisebox{0.00ex} {\footnotesize J}}  \hspace{0.7em} B) =
 \{XA\}  \bigcirc \hspace{-0.9em} {\raisebox{0.00ex} {\footnotesize J}}  \hspace{0.7em} B
\end{align*}

\begin{align*}
(A  \bigcirc \hspace{-0.9em} {\raisebox{0.00ex} {\footnotesize J}}  \hspace{0.7em} B) X =
  A \bigcirc \hspace{-0.9em} {\raisebox{0.00ex} {\footnotesize J}}  \hspace{0.7em} \{BX\}
\end{align*}

Similar rule applies for general compatible tensor multiplication. \\
\item \textbf{Associative :}\\
\begin{align*}
A \bigcirc \hspace{-0.9em} {\raisebox{0.00ex} {\footnotesize J}}  \hspace{0.7em}[B \bigcirc \hspace{-0.9em} {\raisebox{0.00ex} {\footnotesize J}}  \hspace{0.7em}C] =
[A \bigcirc \hspace{-0.9em} {\raisebox{0.00ex} {\footnotesize J}}  \hspace{0.7em}B]\bigcirc \hspace{-0.9em} {\raisebox{0.00ex} {\footnotesize J}}  \hspace{0.7em}C
\end{align*}

\item \textbf{Derivative rule :} \\
\begin{align*}
\dfrac{\partial}{\partial x_i} \{ A(x_1x_2...x_n) \bigcirc \hspace{-0.9em} {\raisebox{0.00ex} {\footnotesize J}}  \hspace{0.7em} B(x_1x_2...x_n)\} =
\dfrac{\partial A}{\partial x_i}  \bigcirc \hspace{-0.9em} {\raisebox{0.00ex} {\footnotesize J}}  \hspace{0.7em} B +
  A \bigcirc \hspace{-0.9em} {\raisebox{0.00ex} {\footnotesize J}}  \hspace{0.7em} \frac{\partial}{\partial x_i} B
\end{align*}
\end{itemize}
\begin{itemize}
\item \textbf{Transpose :}\\
\begin{align*}
C_{ijk}^T=C_{kji}\\
\end{align*}
\begin{align*}
[A \bigcirc \hspace{-0.9em} {\raisebox{0.00ex} {\footnotesize J}}  \hspace{0.7em}B]^T =
[B^T \bigcirc \hspace{-0.9em} {\raisebox{0.00ex} {\footnotesize J}}  \hspace{0.7em}A^T]
\end{align*}

For symmetric A,B             \\
\begin{align*}
[A \bigcirc \hspace{-0.9em} {\raisebox{0.00ex} {\footnotesize J}}  \hspace{0.7em}B] =
[B \bigcirc \hspace{-0.9em} {\raisebox{0.00ex} {\footnotesize J}}  \hspace{0.7em}A]
\end{align*}

\item \textbf{Transpose of triple product :} \\
\begin{align*}
[A \bigcirc \hspace{-0.9em} {\raisebox{0.00ex} {\footnotesize J}} \hspace{-0.5em} \raisebox{-1.50ex} {\tiny p} \hspace{0.7em}B \bigcirc \hspace{-0.9em} {\raisebox{0.00ex} {\footnotesize J}} \hspace{-0.5em} \raisebox{-1.50ex} {\tiny q} \hspace{0.7em} C]^T =
C^T \bigcirc \hspace{-0.9em} {\raisebox{0.00ex} {\footnotesize J}} \hspace{-0.5em} \raisebox{-1.50ex} {\tiny q} \hspace{0.7em} B^T \bigcirc \hspace{-0.9em} {\raisebox{0.00ex} {\footnotesize J}} \hspace{-0.5em} \raisebox{-1.50ex} {\tiny p} \hspace{0.5em}A^T
\end{align*}

Note the reversal of p, q \\
\end{itemize}
\subsection {\textbf{Application of J-product  in the present context}}
Derivative of Edge Matrix w.r.t. co-ordinates of interior point  can be expressed in terms of
J-product , which reduces derivative to an algebraic operation :

\begin{align*}
\frac{\partial A}{\partial w_p} = \frac{z}{2} \left\lbrace I
\bigcirc \hspace{-0.9em} {\raisebox{0.00ex} {\footnotesize J}} \hspace{-0.5em} \raisebox{-1.50ex} {\tiny q}
\hspace{0.7em}A + A
\bigcirc \hspace{-0.9em} {\raisebox{0.00ex} {\footnotesize J}} \hspace{-0.5em} \raisebox{-1.50ex} {\tiny q}
\hspace{0.7em}I \right\rbrace
\end{align*}

Recurrence relation for derivative of kth  power of the Edge Matrix in terms of J-products
of lower powers :\\
\begin{align*}
\frac{\partial A^k}{\partial w_p} = \frac{z}{2} sym \left\lbrace \sum_{i=1}^k A^i \bigcirc \hspace{-0.9em} {\raisebox{0.00ex} {\footnotesize J}} \hspace{-0.5em} \raisebox{-1.50ex} {\tiny q} \hspace{0.7em}A^{k-i} \right\rbrace
\end{align*}

Recurrence relation for kth derivative of  the Edge Matrix in terms of derivatives of lower order e.g. second derivative in terms of first derivative :
\\
%\DeclareMathSizes{13.82}{12.4}{10}{7}
\begin{align*}
\frac{\partial^2 A}{\partial w_p \partial w_q}  = \frac{z}{2} \left\lbrace \frac{\partial A}{\partial w_p} \bigcirc \hspace{-0.9em} {\raisebox{0.00ex} {\footnotesize J}} \hspace{-0.5em} \raisebox{-1.50ex} {\tiny q} \hspace{0.7em}I + I
\bigcirc \hspace{-0.9em} {\raisebox{0.00ex} {\footnotesize J}} \hspace{-0.3em} \raisebox{-1.50ex} {\tiny q}
\hspace{0.7em}\frac{\partial A}{\partial w_p}\right\rbrace
\end{align*}

\subsection{\textbf{Solution of recurrence relations for higher derivatives }}
\begin{eqnarray}
\frac{\partial A^k}{\partial w_p}  &  = & \frac{z}{2} \left\lbrace
W^{k, 1}_{rank(\boldsymbol{\alpha}) \cdot } A^{\alpha_1}
\bigcirc \hspace{-0.9em} {\raisebox{0.00ex} {\footnotesize J}}
\hspace{-0.4em} \raisebox{-1.50ex} {\tiny p}
\hspace{0.7em}A^{\alpha_2} \right\rbrace \\
\frac{\partial^2 A^k}{\partial w_p \partial w_q} & = & \frac{z^2}{4}
\cdot \frac{1}{2!}  \left\lbrace W^{k, 2}_{rank(\boldsymbol{\alpha})}
\left[ A^{\alpha_1}
\bigcirc \hspace{-0.9em} {\raisebox{0.00ex} {\footnotesize J}}
\hspace{-0.4em} \raisebox{-1.50ex} {\tiny p}
\hspace{0.7em}A^{\alpha_2} \bigcirc \hspace{-0.9em} {\raisebox{0.00ex} {\footnotesize J}}
\hspace{-0.4em} \raisebox{-1.50ex} {\tiny q}
\hspace{0.5em} A^{\alpha_3}
  \right. \right. \nonumber \\
&  & ~~~~~~~~~~~~~~~~~~~~~~ \left. \left. +  A^{\alpha_1} \bigcirc \hspace{-0.9em} {\raisebox{0.00ex} {\footnotesize J}}
\hspace{-0.4em} \raisebox{-1.50ex} {\tiny q}
\hspace{0.5em}  A^{\alpha_2}
  \bigcirc \hspace{-0.9em} {\raisebox{0.00ex} {\footnotesize J}}
  \hspace{-0.4em} \raisebox{-1.50ex} {\tiny p}
  \hspace{0.5em} A^{\alpha_3}\right]\right\rbrace
\end{eqnarray}
where
\begin{eqnarray}
\boldsymbol{\alpha} & = & (\alpha_1, \alpha_2) ~~~~~ \alpha_1 + \alpha_2 = k, \alpha_i \ge 0 ~~~~~~~~~~~~~\,\text{for (1)}\nonumber\\
\boldsymbol{\alpha} & = & (\alpha_1, \alpha_2, \alpha_3) ~~~~~ \alpha_1 + \alpha_2 + \alpha_3 = k, \alpha_i \ge 0 ~~\text{for (2)}\nonumber \\
rank(\boldsymbol{\alpha}) & = & \text{number of non-zero ~}\alpha_i's\nonumber\\
W^{k, m}_{rank({\boldsymbol{\alpha}})} & = & \text{fixed weights depending on rank}\nonumber
\end{eqnarray}

Similar formulas hold for higher derivatives.

\section{\textbf{Approximate Evaluation on the Standard Model}}
\begin{itemize}
\item We want to illustrate how a function of distributive type over
  the set of all slack variable can be approximated by superposition
  of exponentials, for various $\psi(s).$
\item For simplicity of exposition consider
\begin{align*}
\psi(s)=\dfrac{1}{s}, \ s > 0 \\
\dfrac{1}{s}=\int_0 ^\infty e^{-sx} dx
\end{align*}
\item In the problem of interest,

all slack variable lie in a bounded interval of real axis.
\item Upper bound : trivial since $ -1 \leq w_i \leq 1$. Hence $\vert
  a^Tw \vert < \vert a \vert _1 \leq n$ for closed walks of length at
  most $n$.
\item Lower bound: at each iteration we round the interior point to
  nearest valid $\pm 1$ solution by a simple method, until optimal
  solution is identified.
\item Hence the minimum value of the slack variable remains bounded
  away from zero.
\end{itemize}

Consider the following nested regions

\[
R_i = \{ x \in \mathbb{R} \mid e^{-i} \leq x \leq e^{j_{max}}\}
~~~~\text{for}~ i = 1, \ldots, L
\]
so that
\[
R_1 \subset R_2 \subset R_3 \cdots \subset R_L
\]

Note: $L$ is not known in advance; we increment it dynamically based
on the number of iterations so far.

Projective invariance of the algorithm implies invariance
w.r.t. simple uniform scaling transformations $T_k: s \rightarrow
e^ks, ~~ k \in \mathbb{Z}$.

\vspace{5pt} To get an efficient approximation exploiting the scale
invariance, we make further substitution
\[
x = e^{-at}
\]

in the integral for $\frac{1}{s}$.

\[
\frac{1}{s} = \int_{- \infty}^{\infty} a\cdot e^{-at} \cdot e^{-
  \left[e^{-at}.s\right]} dt
\]

Approximate the integral by sum, with integer nodes in a bounded range
$\left[-m, M\right]$, $m, M \in \mathbb{N}$.

\[
\frac{1}{s} \sim \sum_{i = -m}^{i = M} a\cdot e^{-ia} \cdot e^{-
  \left[e^{-ia}.s\right]}
\]

Let $\lambda_i = e^{-ia}$

\[
\frac{1}{s} \sim \sum_{i = -m}^{i = M} a \cdot \lambda_i \cdot e^{-
  \lambda_i s} ~~~~~~~~~~~~~~~~~~~(*)
\]

As an example, suppose the region of interest is

\[
e^{-30} < s < e^1
\]

Taking $a = \frac{1}{2}$, $m = \frac{1}{a} = 2$, $M = \frac{30}{a} =
60$, we have 63 terms and this gives the l.h.s. and r.h.s. of $(*)$
which are close when evaluated in standard double precision
(64-bit) arithmetic.

%% \begin{itemize}
%% \item
%% \end{itemize}
\section{\textbf{Superposition of Exponentials}}
\begin{itemize}
\item  By similar techniques approximation to the function $\psi(s)$ of interest is expressed as sum of exponential
\begin{align*}
\psi(s) = \sum_i c_i e^{-\lambda _is}
\end{align*}
\item The main function $\phi$ (either metric or potential) is expressed in terms of $\psi(s)$ in distributive form
\begin{align*}
\phi & = \sum_{s \epsilon S} \psi(s), \\
& = \sum_{s \epsilon S} \sum_i c_i e^{-\lambda _is} \\
& = \sum_i c_i \left\lbrace \sum_{s \epsilon S} e^{-\lambda _is} \right\rbrace\\
& = \sum_i c_i \cdot e^z \cdot tr\{B(z)\} \vert _{z=\lambda _i}
\end{align*}
\item Note that the number of terms is $O(L),$ and evaluation of each
  term is polynomial in $n$.
\end{itemize}
\section{\textbf{An example of application to non-convex, non-linear problem}}
\begin{itemize}
\item \vspace{5pt} Consider set of inequalities of the form
\begin{align*}
(a^T _ix)^2 - (b^T _ix)^2 \leq c, \ c > 0
\end{align*}
The set defined by these inequalities is $(2,2)$-connected. Karmarkar \cite{Kar2}
slack variable $s_i = c - (a^T _ix)^2 + (b^T _ix)^2$
\item With the function $\psi(s): \mathbb{R} \rightarrow \mathbb{R}$
  of interest, we associate another function $\tilde{\psi} :
  \mathbb{R} ^2 \rightarrow \mathbb{R}$ as follows
\begin{align*}
\tilde{\psi} (u,v) = \psi (c-u^2+v^2)
\end{align*}
Let $R$ be the region of interest in the $(u,v)$ plane.

Note: due to the symmetries based on flipping signs of $u$ and $v$, we
can consider only one quarter of the plane.
\end{itemize}

\begin{itemize}
\item Let $\chi _R :$ Characteristic function of $R$. By considering
  2-D inverse Laplace transform of $\tilde{\chi}_R \tilde{\psi}
  (u,v),$ (or by other means) we can obtain finite approximation based
  on exponentials functions in the plane
\begin{align*}
\chi _R \tilde{\psi} \sim \sum_{(\mu _k, \lambda _k)} c_k e^{-\left[
    \mu _ku + \lambda _kv\right] }
\end{align*}
\item Since $u=a^T _ix$ and $v=b^T _ix$ are linear in $x,$ the above
  expression is of the same type as considered before.
\end{itemize}

\section{\textbf{Extension to the satisfiability problem}}

\subsection{\textbf{Introduction to SAT}}
\vspace{10pt}
\begin{itemize}
\item We have boolean variables $x_1, \ldots, x_n$.
\item A literal $d$ is a variable or its complement: $d = x$ or $d = \bar{x}$
\item A clause is an OR of literals. Eg. $C = x \vee y \vee
  \bar{z}$. In 3-SAT, each clause has 3 literals.
\item A satisfiability formula $f = C_1 \wedge C_2 \wedge \ldots
  \wedge C_m$
\item Problem (SAT): Find truth assignments to the variables $x_1,
  \ldots, x_n$ such that $f$ is true, or prove that no such assignment
  exists.
\end{itemize}

\subsection{\textbf{Generation of implied clauses}}
\vspace{5pt}
\begin{itemize}
\item Consider the following pair of clauses
\begin{eqnarray}
C_1 & = & a \vee b \vee \bar{c}\nonumber\\
C_2 & = & c \vee d \vee e\nonumber
%&& \int x^2 = \frac{x^3}{3} + C \label{eq:xdef}
\end{eqnarray}
\item Together they imply the clause $C_3 = a \vee b \vee d \vee
  e$.
\item We refer to this operation as the ``join'' operation on clauses.
\item This operation is fundamental in resolution or elimination
  methods for solving the SAT problem.
\item It is known that the number of such new clauses generated during
  the resolution process grows exponentially with $n$ for almost all
  instances of the problem with parameters in a certain range
  Chv\'atal and Szemer\'edi\cite{Chvatal}.
\end{itemize}

\subsection{\textbf{Continuum based approach}}
\begin{itemize}
\item \vspace{5pt} Associate real variables taking $\pm 1$ values with
  the boolean variables.
\item Embedding the problem in $\mathbb{R}^n$, allow the variables to
  take on values in the interval $\left[-1, 1\right]$, i.e. $-1 \leq
  x_i \leq 1$.
\item For a literal $l$,
\[
v(l) =
\begin{cases}
  x & \mbox{if }l = x\\
  -x & \mbox{if }l = \bar{x}
\end{cases}
\]
\item Each clause corresponds to an inequality
\begin{center}
\begin{tabular}{lllllll}
$C$ & = & $ l_1 \vee l_2 \vee l_3$ & $\longrightarrow$ & $v(l_1) + v(l_2) + v(l_3)$ & $\ge$ & -1\\
\end{tabular}
\end{center}
Explanation: Sum of three $\pm 1$ variables can be -3, -1, 1 or 3 of
which the value -3 is forbidden since it corresponds to all literals
being false.
\end{itemize}

\subsection{\textbf{Interpreting the join operation in the continuum model\\}}
% Let us consider the effect of join operation \\
%\vspace{10pt}

\begin{tabular}{llllllll}
\vspace{0pt} & \underline{Clauses} & & & \underline{Inequalities} & & & \\
\vspace{0pt}& $ a \vee b \vee \bar{c}$ & $(C_1)$ & $\longrightarrow$ & $a + b - c$ & $\ge$ & -1& (1) \\
\vspace{0pt}& $ c \vee d \vee e$ & $(C_2)$ & $\longrightarrow$ & $c + d + e$ & $\ge$ & -1 & (2) \\
\vspace{0pt}Join: & $ a \vee b \vee d \vee e$ & $(C_3)$ & $\longrightarrow$ & $a + b + d + e$ & $\ge$ & -2 & (3) \\
\end{tabular}
\\
\begin{itemize}
\item Observe that the inequality (3) corresponding to the newly
  generated clause is just the sum of (1) and (2).
\item Imposition of constraint (3) does not change the feasible set,
  hence it is superfluous.
\item This is the first reason why the continuum approach is more
  economical than resolution.
\end{itemize}

\subsection{\textbf{A thought experiment}}
\vspace{10pt}
\begin{itemize}
\item Run the resolution algorithm.
\item For each new clause generated, ask if the corresponding
  inequality is superfluous or essential (i.e. not derivable as a
  non-negative combination of the previous inequalities).
\item Unless the previous inequalities have only vectors with all $\pm
  1$ co-ordinates as extreme points, an essential constraint must get
  generated during the course of the algorithm since resolution is a
  complete method.
\item This leads us to the following question: Which sequence of
  clauses when joined, yields an essential constraint, and is such
  that no subsequence has the same property?
\item To understand this, we introduce the concepts of paths, walks,
  cycles, etc. formed by subformulas in a SAT problem.
\end{itemize}

\section{\textbf{Paths, Cycles and Mobius cycles}}
\vspace{5pt}
\subsection{\textbf{(Open) Path}}:
\begin{itemize}
\item \vspace{3pt} Consider a sequence of clauses of the following
  type:

$l_1 \vee l_2 \vee \overline{l_3}$, ~~$\,l_3 \vee l_4 \vee
  \overline{l_5}$, ~~$\,l_5 \vee l_6 \vee \overline{l_7}$, $\ldots,
  l_{2k-1} \vee l_{2k} \vee \overline{l_{2k+1}}$

where each $l_i$  is a literal based on a distinct variable.
\item \vspace{3pt} For any two consecutive clauses in the sequence,
  the last literal of the earlier clause and the first literal of the
  latter clause are complementary.
\item \vspace{3pt} The above sequence yields the following joined
  clause:

\hspace{40pt}$l_1 \vee l_2 \vee l_4 \vee l_6 \vee \ldots \vee l_{2k}
\vee \overline{l_{2k+1}}$

and the corresponding inequality is superfluous.
\end{itemize}
\subsection{\textbf{(Ordinary) Cycle}}:
\begin{itemize}
\item \vspace{3pt} If $l_{2k+1} = l_1$ above, then we can join the two
  ends by normal rules of joining, but the joined clause of the
  sequence is a tautology (always true) due to the presence of $l_1$
  and $\overline{l_1}$.
\end{itemize}

\subsection{\textbf{Mobius cycles}}
\begin{itemize}
\item \vspace{3pt} Consider the above sequence of clauses again:

$l_1 \vee l_2 \vee \overline{l_3}$, $\,l_3 \vee l_4 \vee
  \overline{l_5}$, $\,l_5 \vee l_6 \vee \overline{l_7}$, $\ldots,
  l_{2k-1} \vee l_{2k} \vee \overline{l_{2k+1}}$

\item If $l_{2k+1} = \overline{l_1}$, then the joined clause of the
  sequence contains two copies of $l_1$, but we need only one. Hence,
  the joined clause is equivalent to $l_1 \vee l_2 \vee l_4 \vee l_6
  \vee \ldots \vee l_{2k}$.
\item Corresponding inequality
\begin{center}
$v(l_1) + v(l_2) + v(l_4) + \ldots + v(l_{2k}) \ge -k + 1$
\end{center}
where
\[
v(l) =
\begin{cases}
  x & \mbox{if }l = x\\
  -x & \mbox{if }l = \bar{x}
\end{cases}
\]
\item This is sharper than the inequality implied by the sum of the
  constituent inequalities, namely
\begin{center}
$2 \cdot v(l_1) + v(l_2) + v(l_4) + \ldots + v(l_{2k})  \ge  -k$
\end{center}
\end{itemize}

\vspace{20pt}
\begin{itemize}
\item \vspace{3pt} We call the cycles of the kind referred to above as
  ``mobius cycles''.
\item \vspace{3pt} The mobius cycles are the only mechanism giving
  rise to new sharper inequalities.
\item \vspace{3pt} If there were no mobius cycles, then solving the
  L.P. corresponding to the SAT problem is sufficient to solve the
  latter.
\item \vspace{3pt} In other words, the special classes of SAT formulas
  in which mobius cycles are forbidden as subformulas, can be solved
  in polynomial time.
\end{itemize}

\section{\textbf{Computing the effect of mobius cycles in the continuum model}}
\begin{itemize}
\item The total number of mobius cycles in the original formula can be
  exponentially large, but we can compute their net effect in
  polynomial time.
\item The method is similar to the one used for odd cycles in the
  maximum independent set problem, except that we construct a directed
  graph as explained below.
\item $2n$ nodes corresponding to the literals $x_1, \overline{x_1},
  x_2, \overline{x_2}, \ldots, x_n, \overline{x_n}$.
\item Let $C = l_1 \vee l_2 \vee l_3$ be a clause.
\item Inequality corresponding to this clause $C$
\[
v(l_1) + v(l_2) + v(l_3)  \ge -1
\]
\item Slack variable $s = 1 + v(l_1) + v(l_2) + v(l_3)$

\[
\psi_C(s) = e^{-z \cdot s} = e^{-z} \cdot e^{z \left[v(l_1) + v(l_2) + v(l_3)\right]}
\]

\end{itemize}
\begin{itemize}
\item Corresponding to $C$, put the following 6 edges in the graph.
\begin{center}
\begin{tabular}{l}
$(l_1, \overline{l_2}), (l_1, \overline{l_3}), (l_2,
  \overline{l_3})$\\
$(\overline{l_1}, l_2), (\overline{l_1}, l_3), (\overline{l_2},
  l_3)$\\
\end{tabular}
\end{center}
\item Each of the 6 edges corresponding to clause $C$ is labelled with
  $\psi_C(s)$.
\item Construct a ``clause matrix'' $A(z, \textbf{x})$

\[
A_{i j} =
\begin{cases}
\text{sum of all } \psi_C(s)~\text{of all parallel edges} & \\
\text{between } i~\text{and}~ j & \text{if any}\\
0 & \text{otherwise}\\
\end{cases}
\]
\item Note that the clause matrix is not symmetric
\[
A(i, j) \neq A(j, i)
\]
but it has another kind of symmetry
\[
A(i, j) = A(\bar{j}, \bar{i})
\]

\end{itemize}
%% \begin{itemize}
%% \item
%% \end{itemize}
\section{\textbf{Closed form expression for combined effect of all mobius
  cycles}}
\begin{itemize}
\item As before, we permit self-intersecting, directed, closed walks
  and thereby include the effect of ``useless but harmless''
  inequalities. This leads to closed form expression for the combined
  effect as before.
\item Consider a mobius cycle
$l_1 \vee l_2 \vee \overline{l_3}$, $\,l_3 \vee l_4 \vee
  \overline{l_5}$, $\,l_5 \vee l_6 \vee \overline{l_7}$, $\ldots,
  l_{2k-1} \vee l_{2k} \vee l_1$
\item Let $s$ and $\widetilde{s}$ be the slack variables for the
  implied and sharper inequalities respectively.
\begin{tabular}{lll}
$s$ & = & $k + 2 \cdot v(l_1) + v(l_2) + v(l_4) + \ldots +
v(l_{2k})$\\
$\widetilde{s}$ & = & $k-1 + v(l_1) + v(l_2) + v(l_4) + \ldots +
v(l_{2k})$\\
$s - \widetilde{s}$ &  = & $1 + v(l_1)$\\
$e^{z \cdot (s - \widetilde{s})}$ & = & $e^{z \left[1 + v(l_1)\right]}$\\
$\psi(\widetilde{s})$ & = & $e^{-z \cdot \widetilde{s}} = e^{-z \cdot s}
\cdot e^{z \left[1 + v(l_1)\right]} = e^{z \left[1 + v(l_1)\right]}
\cdot \psi(s)$\\
\end{tabular}
\item The factor $e^{z \left[1 + v(l_1)\right]}$ is associated with a
  special edge from $\overline{l_1}$ to $l_1$. Traversing this edge
  differentiates mobius cycles from ordinary cycles.
\end{itemize}
\begin{itemize}
\item Construct ``mobius completion matrix'' $M_c$ from special edges as follows:
\begin{eqnarray*}
M_C(z, \textbf{x})_{i, \bar{i}} & = & e^{z \left[ 1 - x_i \right]} \nonumber\\
M_C(z, \textbf{x})_{\bar{i}, i} & = & e^{z \left[ 1 + x_i \right]} \nonumber\\
M_C(z, \textbf{x})_{i, j} & = & 0 ~~\text{if } j \neq \bar{i} ~~~~\{\text{Note:}
\bar{\bar{i}} =  i\} \nonumber
\end{eqnarray*}
\item Define walk matrix $B(z, \textbf{x})$ as before, except that we include
  all walks (odd and even) upto length $k_{max}$, where $k_{max} = $
  number of clauses.
\[
B(z, \textbf{x}) = \sum_{k = 2}^{k = k_{max}} \frac{A^k(z, \textbf{x})}{2k}
\]
\item For each walk, there is a ``mirror-image'' walk obtained by
 complementing all the nodes in t
 he walk, hence factor of 2 in the denominator

\item Multiplication by $M_c $ also achieves the transition from implied to sharper inequalities. The combined effect is given by

\[
\phi(z, \textbf{x}) =  tr\{M_c(z, \textbf{x}) \cdot B(z, \textbf{x})\}
\]

\item \vspace{3pt} Lengths of mobius cycles have expontential effect on
resolution. In contrast, we are able to include the effect of mobius
cycles of all lengths in polynomial number of operations.\\

The remaining part of computation is similar to the case of maximum independent set problem. 

\end{itemize}

% -----------------------------------------------------------------------------

\bibliographystyle{amsplain}

% -----------------------------------------------------------------------------

\end{document}